\documentclass[12pt,reqno]{amsart}

\usepackage{amsmath}
\usepackage{amssymb}
\usepackage{amsthm}
\usepackage[english]{babel}

\usepackage{graphics, xcolor}

\usepackage{amsmath}
\usepackage{amssymb}
\usepackage{amsthm}
\newtheorem{theorem}{Theorem}

\newtheorem{lemma}{Lemma}

\newcommand{\beqa}{\begin{eqnarray}}
\newcommand{\beqan}{\begin{eqnarray*}}
\newcommand{\eeqa}{\end{eqnarray}}
\newcommand{\eeqan}{\end{eqnarray*}}
\def\beq#1\eeq{\begin{equation}#1\end{equation}}

\def\P{\mathbf P }

 \def\na{\,\, {\raise.4pt\hbox{$\shortmid$}}{\hskip-2.0pt\to}\, \, }

\def\={\overset{ \text{\rm def} }=}

\def\ffrac{\frac}
\def\4{\kern1pt}

\def\bgr#1{\4\bigr#1}
\def\bgl#1{\bigl#1\4}

\def\gr{}
\def\tc{}%\newcommand{\bl}{}
\newtheorem{remark}{Remark}

\begin{document}
\title[Estimates in Kolmogorov's theorem]
{\gr{On alternative approximating distributions in the  multivariate version of  Kolmogorov's second uniform limit theorem}}

\author[F.~G\"otze]{Friedrich G\"otze}
\author[A.Yu. Zaitsev]{Andrei Yu. Zaitsev}

\email{goetze@math.uni-bielefeld.de}
\address{Fakult\"at f\"ur Mathematik,\newline\indent
Universit\"at Bielefeld, Postfach 100131,\newline\indent D-33501 Bielefeld,
Germany\bigskip}
\email{zaitsev@pdmi.ras.ru}
\address{St.~Petersburg Department of Steklov Mathematical Institute
\newline\indent
Fontanka 27, St.~Petersburg 191023, Russia\newline\indent
and St.Petersburg State University, 7/9 Universitetskaya nab., St. Petersburg,
199034 Russia}

\begin{abstract}
The aim of the present work is to show that recent results of  the authors
on the approximation of distributions of sums of independent summands by
the infinitely divisible laws on convex polyhedra can be shown via an alternative class of approximating infinitely divisible distributions.
We will also generalize the results to the infinite-dimensional case.
\end{abstract}

\keywords {sums of independent random variables, closeness of successive convolutions, convex polyhedra, approximation, inequalities}

\subjclass {Primary 60F05; secondary 60E15, 60G50}

\thanks{The authors were supported by the SFB 1283 and by RFBR-DFG grant 20-51-12004.
The second author was supported by grant RFBR 19-01-00356.}

\maketitle

We would like to show that some of our recent results in \cite{gzz}
may be derived based on an alternative class of infinitely divisible distributions. We will also generalize the results to the infinite-dimensional case.

Let us first introduce some notation. Let $\frak F_d$ denote the set
of probability distributions defined on the Borel $\sigma$-field of
subsets of the Euclidean space~${\mathbf R}^d$. Let $\mathfrak{D}_d\subset\mathfrak{F}_d$ be the set of infinitely divisible distributions. For $F\in\frak F_d$, we
denote the corresponding  distribution functions
by~$F(b)$:
$$
F(b)=F\{(-\infty,b_1]\times\dots\times(-\infty,b_d]\},\qquad b=(b_1,\dots,b_d)\in{\mathbf R}^d.
$$
Let \,${\mathcal L}(\xi)\in\frak F_d$ \,be the
distribution of a $d$-dimensional random vector~$\xi$.
Products and powers of measures are understood in the
convolution sense: $$
GH=G*H,\quad H^m=H^{m*},\quad H^0=E=E_0,
$$
where $E_x$ is the distribution concentrated at a point $x\in
{\mathbf R}^d $.
By~\,$c$ \,we
denote absolute positive constants. Note that constants \,$c$ \,can be different
in different (or even in the same) formulas. If the corresponding constant depends
on, say, $s$, we write $c(s)$.

Kolmogorov \cite{23} posed the problem of estimating the accuracy of infinitely divisible
approximation of distributions of sums of independent random
variables, the distributions of which are concentrated on the short
intervals of length $\tau\le1/2 $ to within a small probability~$p$.
The restriction on the distributions of {summands} is a non-asymptotic analogue of the classical infinitesimality (negligibility) condition for a triangular scheme of independent random
variables. The bound for the rate of approximations may be considered as a quantitative
improvement of the classical Khinchin theorem {for}  the set of infinitely divisible
distributions being  limit laws of the distributions of sums
	in a triangular scheme.
Suppose that the
distributions  $F_{i}\in\frak F_d$, $i=1,\dots ,n$, are represented as mixtures of probability distributions $U_{i},V_{i}\in\frak F_d$:
\begin{equation}
F_{i}=(1-p_{i})U_{i}+p_{i}V_{i},\label{e1}
\end{equation}
where
\begin{equation}
0\le p_{i}\le 1,\quad \int x\,U_{i}\{dx\}=0,\quad U_{i}\left\{
\left\{ x\in \mathbf{R}^{d}:\left\| x\right\| \le \tau \right\}
\right\} =1, \quad\tau\ge0,\label{e2}
\end{equation}
and $V_{i}$ are arbitrary distributions. Denote
\begin{equation}
p=\max_{1\le i\le n}p_{i},\qquad
F=\prod_{i=1}^{n}F_{i}.\label{e3}
\end{equation}

Kolmogorov \cite{23} proved that in the one-dimensional case, for $d=1$, there exists an
infinitely divisible distribution $D$ such that
\begin{equation}
L(F,D) \le c\biggl( p^{1/5}+\tau^{1/2}\ln^{1/4}\frac1{\tau}\biggr),
\label{infdiv2}
\end{equation}
where
\begin{equation} L(F,D)=\inf\big\{\varepsilon:~F(b-\varepsilon)-\varepsilon\leqslant D(b)\leqslant F(b+\varepsilon)+\varepsilon,\quad\hbox{for all } x\in\mathbf{R}\big\},  \label{Lf}\end{equation} is the
L\'evy distance which metrizes the weak convergence of probability distributions.

This proves Khinchin's theorem since weak convergence $F\Rightarrow H$ implies weak convergence $D\Rightarrow H$ as $p\to0$ and $\tau\to0$. The distribution $H$ is
 infinitely divisible as a limit of infinitely divisible distributions $D$. However, Kolmogorov's inequality \eqref{infdiv2} provides good infinitely divisible approximation for fixed small $p$ and $\tau$ even if the distributions of sums involved in the triangular scheme  with $p\to0$ and $\tau\to0$ are not pre-compact.

Conditions \eqref{e1}--\eqref{e3} do not include any moment restrictions since $V_{i}$ are arbitrary distributions.
Note that the statement of Kolmogorov's result \cite{23} is a little bit different, but it is not difficult to verify the equivalence of formulations.
Later, Kolmogorov \cite{k63} returned to this problem and proved the bound
\begin{equation}
L(F,D) \le c\biggl( p^{1/3}+\tau^{1/2}\ln^{1/4}\frac1{\tau}\biggr).
\label{inf5}
\end{equation}
Ibragimov and Presman \cite{ip73} have shown that it is possible to improve this inequality to
\begin{equation}
L(F,D) \le c\biggl( p^{1/3}+\tau^{2/3}\ln\frac1{\tau}\biggr).
\label{fdiv}
\end{equation}
Finally, the optimal bound was derived in Zaitsev and Arak \cite{49}
\begin{equation}
L(F,D) \le c\biggl( p+\tau\ln\frac1{\tau}\biggr).
\label{infdiv4}\end{equation}
The estimate \eqref{infdiv4} was proved by Zaitsev. Moreover, as was shown by Arak, inequality \eqref{infdiv4} is correct in order with respect to $p$ and~$\tau$. As approximating laws, the so-called accompanying infinitely divisible compound Poisson distributions were used. In 1986,  a
joint monograph by Arak and  Zaitsev~\cite{2},
containing a summary of these results, was published in Proceedings of the Steklov Institute of Mathematics.

Zaitsev~\cite{z89} generalized inequality \eqref{infdiv4} to the multidimensional case. He has shown that, for $d\ge1$,
\begin{equation}
L(F,D) \le c(d)\biggl( p+\tau\ln\frac1{\tau}\biggr),
\label{inf3}
\end{equation}
where
\begin{equation} L(F,D)=\inf\big\{\varepsilon:~F(b-\varepsilon\,\mathbf{1})-\varepsilon\leqslant D(b)\leqslant F(b+\varepsilon\,\mathbf{1})+\varepsilon,\quad\hbox{for all } b\in\mathbf{R}^d\big\},  \label{L}\end{equation}
and $\mathbf{1}\in\mathbf{R}^d$ is the vector with all coordinates equal to one.

The multidimensional L\'evy distance between distributions $ G,H\in\frak F_d $ may be also defined as
\[
L (G,H)=\inf \left\{ \lambda :L (G,H,\lambda )\leq \lambda \right\} ,
\]
where
\begin{equation}
L (G,H,\lambda )=\sup_{ b\in\mathbf{R}^d}\max \big\{ G(b)- H (b+\lambda\,\mathbf{1}), H(b)- G(b+\lambda\,\mathbf{1})\big\},\quad \lambda >0.\label{305}
\end{equation}
The Prokhorov distance between distributions $G,H\in\frak F_d$ may be defined as
\[
\pi (G,H)=\inf \left\{ \lambda :\pi (G,H,\lambda )\leq \lambda \right\} ,
\]
where
\[
\pi (G,H,\lambda )=\sup_{X}\max \left\{ G\{X\}-H\{X^{\lambda
}\},H\{X\}-G\{X^{\lambda }\}\right\} ,\quad \lambda >0,
\]
and $X^{\lambda }=\{y\in \mathbf{R}^{d}:\inf\limits_{x\in X}\left\|
x-y\right\| <\lambda \}$ is the $ \lambda $-neighborhood of a Borel set~$ X $ (see, e,g., \cite{Zol83}).

Le Cam~\cite{LC} proposed \tc{to use as a natural  infinitely divisible approximation of $\prod_{i=1}^nF_i$ the accompanying compound Poisson distribution $\prod_{i=1}^n e(F_i)$, where $$
e(H)=e^{-1}\sum_{s=0}^\infty
 \ffrac{H^s}{s!},\quad \mbox{for }H\in\frak F_d.
$$

If\,$F={\mathcal L}(\xi)\in\frak F_d$ \,and ${\mathbf E}\,\|\xi\|^2<\infty$,
then $\Phi(F)\in\frak F_d $ denotes below the Gaussian distribution with the same mean and covariance operator
as~$F$.
\bigskip

The following Theorem~\ref{th1} is the main result of Zaitsev~\cite{z89}.\bigskip

\begin{theorem}\label{th1}
Let conditions \eqref{e1}--\eqref{e3} be satisfied. Denote
\begin{equation}
 D=\prod_{i=1}^{n}\mbox{e}(F_{i}),\label{e4}
\end{equation}
Then, for any $\lambda>0$,
\begin{equation}
L(F,D,\lambda) \le
c(d)\,\Big( p+\exp \Big( -\frac{\lambda
}{c(d)\,\tau }\Big) \Big),
\label{infdiv}
\end{equation}
\begin{equation}
\pi(F,D,\lambda) \le
c(d)\,\Big( p+\exp \Big( -\frac{\lambda
}{c(d)\,\tau }\Big) \Big) +\sum_{i=1}^{n}p_{i}^{2}.
\label{infdiv1}
\end{equation}
Hence,
\begin{equation}
L(F,D) \le
c(d)\,\left( p+\tau(|\ln \tau|+1) \right),
\label{kkk}
\end{equation}
\begin{equation}
\pi(F,D) \le
c(d)\,\left(p+\tau(|\ln \tau|+1) \right) +\sum_{i=1}^{n}p_{i}^{2}.
\label{infdiv3}
\end{equation}
Inequalities \eqref{infdiv}--\eqref{infdiv3} remain true after replacing $D$ by other approximating infinitely divisible distributions
\begin{equation}\label{d1}
 D^*=\Phi\Big(\prod_{i=1}^{n}\big((1-p_{i})U_{i}+p_{i}E\big)\Big)\prod_{i=1}^{n}{e}
 \big((1-p_{i})E+p_{i}V_{i}\big)
\end{equation}
 or
\begin{equation}\label{d2}
 D^{**}=D_0\,\prod_{i=1}^{n}{e}\big((1-p_{i})E+p_{i}V_{i}\big),
\end{equation}
 where $D_0$ is an arbitrary infinitely divisible distribution with spectral measure concentrated on
the ball $\left\{ x\in \mathbf{R}
^{d}:\left\| x\right\| \le \tau \right\} $ and with the same mean and the same covariance operator
as those of the distribution $\prod_{i=1}^{n}\big((1-p_{i})U_{i}+p_{i}E\big)$.
\end{theorem}

\begin{remark}\rm
  Formally, similarly to the case $d=1$, we consider Gaussian laws as infinitely divisible distributions with spectral measures concentrated at zero. Thus, the distribution $D_0$ may have a Gaussian component.
\end{remark}

\begin{remark}\rm
It is easy to see that the distributions $D$ and $D^*$ are particular cases of distribution $D^{**}$ with
$$D_0=\prod_{i=1}^{n}e\big((1-p_{i})U_{i}+p_{i}E\big)\quad\mbox{and}\quad D_0=\Phi\Big(\prod_{i=1}^{n}\big((1-p_{i})U_{i}+p_{i}E\big)\Big)
$$
respectively.
\end{remark}

\begin{remark}\rm
The mean and the covariance operator of distribution $D_0$ \tc{may be not precisely  equal to  those of the distribution $\prod_{i=1}^{n}((1-p_{i})U_{i}+p_{i}E)$ but may be just close to them}. The additional remainder term will come from the estimation of the closeness of Gaussian laws
$\Phi\big(D_0\big)$ and $\Phi\big(\prod_{i=1}^{n}((1-p_{i})U_{i}+p_{i}E)\big)$ (see, e.g., \cite{GNU}).
\end{remark}

Note that the estimation of $L (F,D,\lambda )$ and $\pi (F,D,\lambda )$ for all $\lambda>0$ provides more information on the closeness of distributions $F$ and $D$ than the estimation of $L (F,D)$ and $\pi (F,D)$. For example, inequalities \eqref{kkk} and \eqref{infdiv3} are trivial for $\tau\ge1$ while inequalities \eqref{infdiv} and \eqref{infdiv1} are interesting for any $\tau>0$. Moreover, the information containing in \eqref{infdiv} and \eqref{infdiv1} remains invariant if we multiply the random vectors by a non-zero constant.
However, inequalities \eqref{infdiv} and \eqref{infdiv1} actually can be derived from inequalities \eqref{kkk} and \eqref{infdiv3} by varying normalizing factors (see \cite{z84} for details).

Kolmogorov \cite{23, k63} has obtained actually the bounds for $L (F,D^*,\lambda )$, $\lambda\ge2\tau>0$, in the case $d=1$. Instead of \eqref{infdiv2} and \eqref{inf5}, he has proved inequalities
\begin{equation}
L(F,D^*,\lambda ) \le
c\,\Big( p^{1/5}+\frac\tau\lambda \ln^{1/2}\frac \lambda\tau\Big),
\label{in2}
\end{equation}
and
\begin{equation}
L(F,D^*,\lambda ) \le
c\,\Big( p^{1/3}+\frac\tau\lambda \ln^{1/2}\frac \lambda\tau\Big)
\label{in5}
\end{equation}
respectively. The optimality of inequality \eqref{infdiv} means that the case where $\lambda<2\tau$ is trivial: if $\lambda<2\tau$, then there exists $F$ from \eqref{e3} such that $L(F,D,\lambda )\ge c$, for any $D\in \mathfrak{D}_1$.

\bigskip

The proof of Theorem~\ref{th1} is based on the following Lemmas \ref{l1}--\ref{l4}.

\begin{lemma}[see \cite{Zol}]\label{l1} Let $ F, G, H\in\mathfrak F_d$ be arbitrary distributions. Then
$R(FH, GH)\le R(F,G)$, where $R(\,\cdot\,,\,\cdot\,)$ is any of the distances
$L(\,\cdot\,,\,\cdot\,)$, $\pi(\,\cdot\,,\,\cdot\,)$ or $\rho(\,\cdot\,,\,\cdot\,)$ $(\rho$ is the uniform distance between ditribution functions$)$.
\tc{Moreover, $L(\,\cdot\,,\,\cdot\,)\le\min\big\{\pi(\,\cdot\,,\,\cdot\,),\rho(\,\cdot\,,\,\cdot\,)\big\}$.}
\end{lemma}

\begin{lemma}[Zaitsev \cite{z85}]\label{l2} Let the conditions of Theorem~$\ref{th1}$ be satisfied. Let
\[
G_{i}=(1-p_{i})E+p_{i}V_{i},\quad H_{i}=(1-p_{i})U_{i}+p_{i}E,\quad
i=1,\dots ,n,
\]
and
\[
G=\prod_{i=1}^{n}G_{i},\qquad H=\prod_{i=1}^{n}H_{i}.
\]
Then
\begin{equation}
\pi(F,GH) \le
c(d)\,\big(p+\tau(|\ln \tau|+1) \big) .
\label{in4}
\end{equation}
and, for any $\lambda>0$,
\begin{equation}
\pi(F,GH,\lambda) \le
c(d)\,\Big( p+\exp \Big( -\frac{\lambda
}{c(d)\,\tau }\Big) \Big) .
\label{inf1}
\end{equation}
\end{lemma}

\begin{lemma}[Zaitsev \cite{z83}, \cite{z89}]\label{Th1} Assume that the
 distributions $G_i\in\mathfrak F_d$ are represented as
 \begin{equation}
G_i=(1-p_i)\,E+p_i\,V_i,\quad
i=1,\dots ,n, \label{111h}\end{equation}
where $V_i\in\mathfrak F_d$ are arbitrary distributions, $0\le p_i\le p=\max_jp_j$,
$$ G=\prod_{i=1}^nG_i,\quad D=\prod_{i=1}^ne(G_i).$$
Then
\begin{equation} \rho(G,D)  \le c(d)\, p,  \label{r699}\end{equation}
\end{lemma}

\begin{lemma}[{see \cite[p.~186]{LC}}]\label{l7} Let the conditions of Lemma~$\ref{Th1}$ be satisfied.
Then
\begin{equation} \pi(G,D)\le\rho_{\rm TV}(G,D)  \le \sum_{i=1}^{n}p_{i}^{2},  \label{r677}\end{equation}
where
\[
\rho_{\rm TV} (G,D )=\sup_{X}\; \bigl|\,G\{X\}-D\{X
\}\bigr|,
\]
is the distance in total variation
and the supremum is taken over all Borel sets~$ X \subset\mathbf{R}^{d}$.
\end{lemma}

\begin{lemma}[see Zaitsev \cite{z86} or \cite{z87}]\label{l3} Let the conditions of Theorem~$\ref{th1}$ be satisfied with $p=0$, that is
\begin{equation}
\int x\,F_{i}\{dx\}=0,\quad F_{i}\left\{
\left\{ x\in \mathbf{R}^{d}:\left\| x\right\| \le \tau \right\}
\right\} =1, \quad\tau\ge0.
\end{equation}
Then
\begin{equation}
\pi(F,\Phi(F)) \le
c(d)\,\tau(|\ln \tau|+1)
\label{in5}
\end{equation}
and, for any $\lambda>0$,
\begin{equation}
\pi(F,\Phi(F),\lambda) \le
c(d)\,\exp \Big( -\frac{\lambda
}{c(d)\,\tau }\Big)  .
\label{in1}
\end{equation}
\end{lemma}

\begin{lemma}[see Zaitsev \cite{z86} or \cite{z87}]\label{l4} Let $D$ be
an infinitely divisible distribution with spectral measure concentrated on
the ball $\left\{ x\in \mathbf{R}
^{d}:\left\| x\right\| \le \tau \right\} $. Then
\begin{equation}
\pi(D,\Phi(D)) \le
c(d)\,\tau(|\ln \tau|+1)
\label{infdiv6}
\end{equation}
and, for any $\lambda>0$,
\begin{equation}
\pi(D,\Phi(D),\lambda) \le
c(d)\,\exp \Big( -\frac{\lambda
}{c(d)\,\tau }\Big) .
\label{tn1}
\end{equation}
\end{lemma}

\begin{remark}\rm
  The approximating distributions $D^{**}$ were not included in the statement of \cite[Theorem~1.1]{z89} \tc{ but inequalities \eqref{infdiv}--\eqref{infdiv3} are obviously extended to them in view of Lemmas \ref{l1} and~\ref{l4}.}
\end{remark}

Inequality~\eqref{infdiv} is equivalent to the validity of inequality
\begin{equation} \max\Big\{F\{P\}-D\{P_\lambda\},\;D\{P\}-F\{P_\lambda\}\Big\}
\le c(d)\,\Big( p+\exp \Big( -\frac{\lambda
}{c(d)\,\tau }\Big) \Big) \label{g97}\end{equation}
for any $\lambda>0$ and for all sets $P$ and $P_\lambda$ of the form
\begin{equation} P=\big\{x\in\mathbf R^d:\langle x,e_j\rangle\le b_j, \ j=1,\ldots, d\big\}, \label{g9}\end{equation}
and \begin{equation} P_\lambda=\big\{x\in\mathbf R^d:\langle x,e_j\rangle\le b_j+\lambda, \ j=1,\ldots, d\big\}, \label{g98} \end{equation}
where $e_j\in \mathbf R^d$ are the vectors of the standard Euclidean basis, $-\infty< b_j\le \infty$, $j=1,\ldots, d$.

\tc{It is easy to see that $P^\lambda\subset P_\lambda\subset P^{\lambda\sqrt d+\varepsilon}$ for $\varepsilon>0$. Therefore, \eqref{g97} is equivalent to the validity of inequality
\begin{equation} \max\Big\{F\{P\}-D\{P^\lambda\},\;D\{P\}-F\{P^\lambda\}\Big\}
\le c(d)\,\Big( p+\exp \Big( -\frac{\lambda
}{c(d)\,\tau }\Big) \Big) ,\label{g997}\end{equation}
for any $\lambda>0$. In the paper of G\"otze,  Zaitsev and Zaporozhets~\cite{gzz}, it was shown that inequality \eqref{g997}  is valid for convex polyhedra~ $P\in\mathcal P_m$ (see\ \eqref{wq}) with $c(d)$
replaced by $c(m)$ depending only on $m$, the number of half-spaces involved in the definition of a polyhedron~$P$.}

In Theorem \ref{th8} of the present paper, we show that the same statement remains true after replacing $D$ by approximating distributions
$D^*$ and $D^{**}$ from \eqref{d1} and \eqref{d2}. Thus, there is a freedom in the choice of $D_0$ in the definition of appproximating distribution $D^{**}$. The only restriction is that $D_0$ must be an infinitely divisible distribution with spectral measure concentrated on
the ball $\left\{ x\in \mathbf{R}
^{d}:\left\| x\right\| \le \tau \right\} $ and with the same mean and the same covariance operator
as those of the distribution $\prod_{i=1}^{n}\big((1-p_{i})U_{i}+p_{i}E\big)$.
The definition \eqref{L} of a multivariate version of the L\'evy distance is actually not quite natural since the collection of sets $P$ of the form \eqref{g9} is not invariant with respect to rotation while the conditions of Theorem~\ref{th1} are invariant. Therefore, inequality \eqref{g97} remains true after replacing the sets  $P$ and $P_\lambda$ by $\mathbb UP$ and $\mathbb UP_\lambda$, where $\mathbb U$ is a unitary linear operator.
A question is: how to define a multivariate version of the L\'evy distance which can be used in more adequate bounds under the conditions of  Theorem~\ref{th1}?

 In the present paper, we give similar bounds for comparing quantities defined via multivariate polyhedra in G\"otze,  Zaitsev and Zaporozhets~\cite{gzz}. .

For $m\in\mathbf N$ we denote by $\mathcal P_m$ the collection of sets $P \subset \mathbf R^d$
representable in the form
\begin{equation}
\label{wq}
P=\big\{x\in\mathbf R^d: \langle x,t_j\rangle\le b_j, \ j=1,\ldots, m\big\},
\end{equation}
where $t_j\in \mathbf R^d$ are the vectors satisfying $\|t_j\|=1$, $-\infty<  b_j\le \infty$, \ $j=1,\ldots, m$.
The elements of the set $\mathcal P_m $ will be called
{\it convex polyhedra}. They can be  unbounded sets.
For $P\in\mathcal P_m$ defined in \eqref{wq} and $\lambda \ge0$, we denote
\begin{equation} P_\lambda=\big\{x\in\mathbf R^d:\langle x,t_j\rangle\le b_j+\lambda, \ j=1,\ldots, m\big\}. \label{g9877}\end{equation}
By definition, $ P_\lambda$ is the intersection of closed $\lambda$-neighborhoods of half-spaces $\big\{x\in\mathbf R^d:\langle x,t_j\rangle\le b_j\big\}$, $j=1,\ldots, m$.
Clearly, $P^\lambda\subset P_\lambda$. However, $ P_\lambda$ may be essentially larger than $P^\lambda$. For example, it is the case for $m=2$, if $1-\varepsilon<\left|\langle t_1,t_2\rangle\right|<1$ with a small $\varepsilon>0$. In this case the hyperplanes $\big\{x\in\mathbf R^d: \langle x,t_j\rangleß=b_j \big\}$, $j´=1,2$, are almost parallel and the point $x_0$ such that $ \langle x_0,t_j\rangleß=b_j+\lambda$, $j´=1,2$, belongs to $P_\lambda$ and is far from the set $P$. In the proof of Theorem~\ref{th8}  below we will need, however, the inclusion  $P_\lambda\subset P^{c\lambda}$. For this purpose, we will modify the definition of $P_\lambda$.
It is evident that we can rewrite the definition of the polyhedron $P$ adding in it extra restrictions
\begin{equation}
\label{89}
P=\big\{x\in\mathbf R^d: \langle x,t_j\rangle\le b_j, \ j=1,\ldots, m_0\big\},
\end{equation}intersecting $P$ with half-spaces
$H(t_J, b_j)=\big\{x\in\mathbf R^d: \langle x,t_j\rangle\le b_j\big\}$,  $j=m+1,\ldots, m_0$. It will be the same polyhedron if $P\subset H(t_J, b_j)$. for all $j=m+1,\ldots, m_0$.

.

Similarly to \eqref{g9877}, we denote \begin{equation} P_\lambda=\big\{x\in\mathbf R^d:\langle x,t_j\rangle\le b_j+\lambda, \ j=1,\ldots, m_0\big\}. \label{0000}\end{equation}This is the same notation, but here we considered $P$ as an element of $\mathcal P_{m_0}$.The polyhedron $ P_\lambda$ is again the intersection of closed $\lambda$-neighborhoods of half-spaces $\big\{x\in\mathbf R^d:\langle x,t_j\rangle\le b_j\big\}$, $j=1,\ldots, m_0$.  The only difference is that in \eqref{0000} we have more intersecting half-spaces. We choose these half-spaces with $j=m+1,\ldots, m_0$ so that we "cut" points of $P_\lambda$ which were  far from the set $P$ (see Lemma \ref{l7} below which is proved in G\"otze,  Zaitsev and Zaporozhets~\cite{gzz}).

\begin{lemma}\label{l7}Fix some $m\in\mathbf N$ and $\varepsilon>0$.
Let the polyhedron $P\in\mathcal P_{m}$ be defined in~\eqref{wq}. Then there exist  a
 $c_{m,\varepsilon}$ depending on  $m$ and $\varepsilon$ only, $m_0\in\mathbf N$,
  $m_0\leq c_{m,\varepsilon}$, $t_j\in \mathbf R^d$ with $\|t_j\|=1$,  and $b_j\in \mathbf R$, $j=m+1,\dots, m_0$, such that,
 for any $\lambda>0$,
\begin{align}\label{2114}
P_\lambda:=\big\{x\in\mathbf R^d: \langle x,t_j\rangle\le b_j+\lambda, \ j=1,\ldots, m_0\big\}\subset P^{(1+\varepsilon)\lambda}.
\end{align}
\end{lemma}
The statement of Lemma \ref{l7} is almost evident for $d=2$ and $d=3$.

{Following ~\cite{gzz},} define, for $m\in\mathbf N$,  $ G,H\in\mathfrak F_d $,
\begin{equation}
L_m (G,H)=\inf \left\{ \lambda :L_m (G,H,\lambda )\leq \lambda \right\} ,
\end{equation}
where
\[
L_m (G,H,\lambda )=\sup_{P\in\mathcal P _m}\max \left\{ G\{P\}-H\{P_{\lambda
}\},H\{P\}-G\{P_{\lambda }\}\right\} ,\quad \lambda >0.
\]
Define also
\[
\pi_m (G,H)=\inf \left\{ \lambda :\pi_m (G,H,\lambda )\leq \lambda \right\} ,
\]
where
\[
\pi_m (G,H,\lambda )=\sup_{P\in\mathcal P _m}\max \left\{ G\{P\}-H\{P^{\lambda
}\},H\{P\}-G\{P^{\lambda }\}\right\} ,\quad \lambda >0.
\]
\bigskip

\begin{remark}\rm
With a fixed $ m $, it is easy to verify that $L_m (\,\cdot\,,\,\cdot\, )$ is a distance in the space $\mathfrak F_d$,
An open question is to check that for $\pi_m (\,\cdot\,,\,\cdot\, )$. For $m>1$,
it is problematic to prove or disprove the fulfillment of the triangle inequality.
The difficulty is that the ${\lambda }$-neighborhood $P^{\lambda }$ of a convex polyhedron $P$ unlike $P_{\lambda }$ generally speaking is not a convex polyhedron. It is also clear that the distance $L_1 (\,\cdot\,,\,\cdot\, )=\pi_1 (\,\cdot\,,\,\cdot\, )$ metrizes weak convergence.
 \end{remark}

The following Theorems~\ref{th2}--\ref{th4} are the main results of this paper.\bigskip

\begin{theorem}\label{th2}  Let the conditions of Theorem~$\ref{th1}$ be satisfied. Then, for any
$m\in\mathbf N$,
\begin{equation}
L_m(F,D) \le
c(m)\,\big(p+\tau(|\ln \tau|+1)\big),
\label{div2}
\end{equation}
and
\begin{equation}
L_m(F,D,\lambda) \le
c(m)\,\Big( p+\exp \Big( -\frac{\lambda
}{c(m)\,\tau }\Big) \Big)  ,\quad \lambda >0.
\label{div}
\end{equation}
Inequalities \eqref{div2} and \eqref{div} remain true after replacing $D$ by approximating distributions \eqref{d1} and \eqref{d2}.
\end{theorem}

\begin{theorem}\label{th8}   Let the conditions of Theorem~$\ref{th1}$ be satisfied. Then, for any
$m\in\mathbf N$,
\begin{equation}
\pi_m(F,D) \le
c(m)\,\big(p+\tau(|\ln \tau|+1)\big),
\label{adiv2}
\end{equation}
and
\begin{equation}
\pi_m(F,D,\lambda) \le
c(m)\,\Big( p+\exp \Big( -\frac{\lambda
}{c(m)\,\tau }\Big) \Big)  ,\quad \lambda >0.
\label{adiv}
\end{equation}
Inequalities \eqref{adiv2} and \eqref{adiv} remain true after replacing $D$ by approximating distributions $D^*$ and $D^{**}$ from \eqref{d1} and \eqref{d2}.
\end{theorem}

Thus, the statement of Theorem~\ref{th1} is generalized, since Theorems~\ref{th2} and~\ref{th8} deal with the values of distributions on convex polyhedra \eqref{wq} whereas Theorem~\ref{th1} corresponds to the sets~\eqref{g9}. Note also that in Theorem~\ref{th1} the constants depend on the dimension~$d$, while in Theorems~\ref{th2} and~\ref{th8} the constants depend only on~$m$ involved in the definition of polyhedra \eqref{wq}. Note that in G\"otze,  Zaitsev and Zaporozhets~\cite{gzz} we have proved
Theorems~\ref{th2} and~\ref{th8} for approximating distributions $D$ only.

\tc{The proof of Theorem~\ref{th2} is  based on applying the $m$-variate version of Theorem~\ref{th1}.
Indeed, the $m$-variate vectors with coordinates $\langle \xi,t_j\rangle$, $\langle \eta,t_j\rangle$, $t_j\in\mathbf R^d$, $\|t_j\|=1$,
$j=1,\ldots, m$, satisfy actually the same $m$-dimensional conditions as the random vectors $\xi,\eta\in\mathbf R^d$ with compared $d$-dimensional distributions $F$ and $D$ from Theorem \ref{th1}.
Let $\mathbb{A}:\mathbf{R}^d\to\mathbf{R}^m$ be the linear operator mapping $x\in \mathbf{R}^d$ to the vector with coordinates $\langle x,t_j\rangle$,
$j=1,\ldots, m$. The vectors $\mathbb{A}\xi,\mathbb{A}\eta$ satisfy the conditions of $m$-variate version of
Theorem~\ref{th1} with replacing $\tau$ by $\tau\sqrt m$. This follows from inequality $\|\mathbb{A}\|\le\sqrt m$.
Thus, roughly speaking, from the known estimates of the distance $L$ in space $\mathbf R^m$ we derive estimates of the distance $L_m$ in  $\mathbf R^d$. Theorem~\ref{th8} will be derived from Theorem~\ref{th2} with the help of Lemma~\ref{l7}.}
\medskip

It is not difficult to understand that the conditions of Theorems~\ref {th2} and~\ref {th8} are meaningful even for $d=\infty$, that is, for distributions in the Hilbert space $\mathbf R^\infty=\mathbf{H}$. The definitions of $L_m (\,\cdot\,,\,\cdot\, )$ and $\pi_m (\,\cdot\,,\,\cdot\, )$
are applicable to such distributions without changes.
\medskip

\begin{theorem}\label{th4}The statements of   Theorems~$\ref {th2}$ and~$\ref {th8}$  remain true for $d=\infty$.
\end{theorem}

Theorem~\ref {th4} can be considered as an adequate infinite-dimensional version of Kolmogorov's second uniform limit theorem. Recall that inequality \eqref{infdiv4} (and hence inequalities~\eqref{div2} and~\eqref{adiv2}) are correct in order with respect to parameters $p$ and~$\tau$.\medskip

It is possible, for example, to use Theorem \ref {th4} for comparing the distributions of random polygonal lines constructed via partial sums of independent random variables with distributions of accompanying processes with independent increments.
\medskip

\begin{remark}\rm
In the authors' papers~\cite {GZ17}, some bounds for the distance $\rho (\,\cdot\,,\,\cdot\, )$ were transferred to the distance $\rho_m (\,\cdot\,,\,\cdot\, )$, $m\in\mathbf N$, defined by equality
$$ \rho_m (F,G)=\sup_{P\in\mathcal P_m} \bigl|F\{P\}-G\{P\}\bigr|, \quad F, G\in\mathfrak{F}_d. $$
In particular, this was done for the assertion of Lemma~\ref {Th1} of the present paper.
This and other results of~\cite{GZ17} are also transferred to the infinite-dimensional case $d=\infty$.
The authors have devoted  a recent publication~\cite {GZ21} to this topic.
\end{remark}
\bigskip\emph{}

 \noindent {\it Proof of Theorem\/ $\ref
 {th2}$. }Fix some  polyhedron $P\in\mathcal P_m$:
\begin{align*}
    P=\big\{x\in\mathbf R^d: \langle x,t_j\rangle\le b_j, \ j=1,\ldots, m\big\}.
\end{align*}where $t_j\in\mathbf R^d$, $\|t_j\|=1$, $b_j\in\mathbf R$,
$j=1,\ldots, m$.
Let $\mathbb A:\mathbf R^d\to\mathbf R^m$ be a linear operator mapping as
\begin{align*}
    x\mapsto y=\big(\langle x,t_1\rangle,\dots,\langle x,t_m\rangle\big).
\end{align*}
Let $e_1,\dots, e_m$ are the vectors of the standard Euclidean basis in $\mathbf R^m$. Consider the polyhedron $\widetilde{P}\subset\mathbf R^m$ belonging to the class $\mathcal P_m^*$ of sets of the form~\eqref{g9} with $d=m$ and defined as
\begin{align*}
    \widetilde{P} = \big\{y\in\mathbf R^m:\langle y,e_j\rangle\le b_j, \ j=1,\ldots, m\big\}.
\end{align*}
Since
\begin{align*}
    \langle x,t_j\rangle=\langle x,\mathbb A^* e_j\rangle=\langle \mathbb Ax,e_j\rangle,
\end{align*}with adjoint operator $\mathbb A^*:\mathbf R^m\to\mathbf R^d$,
this implies that, for any random vector $\xi\in\mathbf R^d$, we have
\begin{align*}
    \P[\xi\in P]=  \P[\mathbb A\xi\in \widetilde{P}]\quad\text{and}\quad\P[\xi\in P_{\lambda}]=  \P[\mathbb A\xi\in \widetilde{P}_{\lambda}].
\end{align*}where
\begin{align*}
    P_{\lambda}=\big\{x\in\mathbf R^d: \langle x,t_j\rangle\le b_j+{\lambda}, \ j=1,\ldots, m\big\}.
\end{align*}
Hence, for any random vectors $\xi,\xi'\in\mathbf R^d$ we have
\begin{align}
&\max\{\mathbf{P}[\xi\in P]-\mathbf{P}[\xi'\in P_\lambda],\, \mathbf{P}[\xi'\in P]-\mathbf{P}[\xi\in P_\lambda]\}
\nonumber
\\
&\qquad =\max\bigl\{\mathbf{P}[\mathbb A\xi\in \widetilde{P}]-\mathbf{P}[\mathbb A\xi'\in \widetilde{P}_\lambda],\, \mathbf{P}[\mathbb A\xi'\in \widetilde{P}]-\mathbf{P}[\mathbb A\xi\in \widetilde{P}_\lambda]\bigr\}
\nonumber
\\
&\qquad\le L(\mathcal L(\mathbb A\xi),\mathcal L\bigl(\mathbb A\xi'),{\lambda}\bigr),
\label{14}
\end{align}
where in the last step we used~\eqref{305}.

 The distributions of $m$-variate vectors with coordinates $\langle \xi,t_j\rangle$, $\langle \eta,t_j\rangle$, $t_j\in\mathbf R^d$, $\|t_j\|=1$,
$j=1,\ldots, m$, actually satisfy the same $m$-dimensional conditions as the distributions of random vectors $\xi,\eta\in\mathbf R^d$ with compared $d$-dimensional distributions $F$ and $D^{**}$, $D^{*}$ or $D$ from Theorem \ref{th1}.
 Indeed,
 let  $\alpha_i\in\mathbf{R}$, $X_i, Y_i\in\mathbf{R}^d$, $i=1,\ldots,n$, be independent random variables and vectors such that
\begin{equation}
\P[\alpha_i=1]=1-\P[\alpha_i=0]=p_i,\quad
\mathcal L(X_i)=U_i. \quad\mathcal L(Y_i)=V_i.\quad i=1,\dots, n.\label{e187}
\end{equation}Let\begin{equation}
\quad \xi_i=(1-\alpha_i)X_i+\alpha_i Y_i,\quad i=1,\dots, n.\label{e1877}
\end{equation}Then\begin{equation}
\mathcal L(\xi_i)=F_i=(1-p_{i})U_{i}+p_{i}V_{i},\label{e1877}
\end{equation}\begin{equation}
\mathcal L(\mathbb A\xi_i)=F_i^{(\mathbb A)}=(1-p_{i})U_{i}^{(\mathbb A)}+p_{i}V_{i}^{(\mathbb A)},\quad i=1,\dots, n.\label{e1877}
\end{equation}Here and below for $W=\mathcal L(\xi)\in \mathfrak F_d$ we write $W^{(\mathbb A)}=\mathcal L(\mathbb A\xi)\in \mathfrak F_m$. If $W$ is an  infinitely divisible distribution with spectral measure concentrated on
the ball
 $$
\{ x\in \mathbf{R}^d\colon \|x\| \le \tau\},
$$
then  $W^{(\mathbb A)}$ is an  infinitely divisible distribution with spectral measure concentrated on
the ball $\left\{ x\in \mathbf{R}^m:\left\| x\right\| \le \tau\sqrt m\right\}$. It suffices to verify that for $W=e(\lambda E_{\tau e})$, $\lambda\ge0$, $e\in \mathbf R^d$,  $\left\| e\right\|=1$. Then $W^{(\mathbb A)}=e(\lambda E_{\tau\mathbb A e})$. It is easy to see that $(e(W))^{(\mathbb A)}= e(W^{(\mathbb A)})$.  It remains to note that $\|\mathbb{A}\|\le\sqrt m$.
Similarly, using \eqref{e2}, we see that \begin{equation}
\int x\,U_{i}^{(\mathbb A)}\{dx\}=0,\quad U_{i}^{(\mathbb A)}\left\{
\left\{ x\in \mathbf{R}^{d}:\left\| x\right\| \le \tau\sqrt m \right\}
\right\} =1,\quad i=1,\dots, n.\label{e25}
\end{equation}
 If the vectors $ \xi, \xi'$ have the same covariance operators, then the covariance operators of the vectors $ \mathbb{A}\xi, \mathbb{A}\xi'$ coinside too. Thus, the distributions $F^{(\mathbb{A})}, D^{(\mathbb{A})}$, $(D^{**})^{(\mathbb{A})}$ satisfy the conditions of $m$-variate version of
Theorem~\ref{th1} imposed on $F$,  $D$,  $D^{**}$ when replacing $\tau$ by $\tau\sqrt {m}$.  Applying
Theorem~\ref{th1}, we obtain, for any $\lambda>0$,
\begin{equation} L(F^{(\mathbb{A})},(D^{**})^{(\mathbb{A})},\lambda)
\le  c(m)\, \Big( p+\exp \Big( -\frac{\lambda
}{c(m)\,\tau }\Big) \Big). \label{999}\end{equation} Using \eqref{14} and \eqref{999}, we come to inequality
\begin{equation}\label{ppp}
L_m(F,D^{**},\lambda) \le
c(m)\,\Big( p+\exp \Big( -\frac{\lambda
}{c(m)\,\tau }\Big) \Big)  ,\quad \lambda >0.
\end{equation}Recall that distributions $D$ and $D^*$ are particular cases of the distribution $D^{**}$.
The  second inequality of Theorem~\ref{th2} follows now from \eqref{ppp}. The first inequality follows from the second one by standard arguments.  Theorem~\ref{th2} is proved.
\medskip

\noindent {\it Proof of Theorem\/ $\ref
 {th8}$.}
Fix some  polyhedron $P\in\mathcal P_m$:
\begin{align*}
    P=\big\{x\in\mathbf R^d: \langle x,t_j\rangle\le b_j, \ j=1,\ldots, m\big\}.
\end{align*}
It follows from Lemma~\ref{l7} that it is possible to represent $P$ in the form
\begin{align*}
    P=\big\{x\in\mathbf R^d: \langle x,t_j\rangle\le b_j, \ j=1,\ldots, m_0\big\}.
\end{align*}
such that
\begin{align*}
    P_{\lambda/2}\subset P^{\lambda}\quad\text{and}\quad m_0\leq N_m\in\mathbf N,
\end{align*}
where
\begin{equation}
P_{\lambda/2}=\biggl\{x\in\mathbf R^d\colon \langle x,t_j\rangle\le b_j+\frac{\lambda}2, \, j=1,\dots, m_0\biggr\}
\end{equation}
and
the constant $N_m$ depends on $m$ only. Thus for any random vectors $\xi,\xi'\in\mathbf R^d$ we have
\gr{\begin{align*}
&\max\bigl\{\mathbf{P}[\xi\in P]-\mathbf{P}[\xi'\in P^\lambda],\, \mathbf{P}[\xi'\in P]-\mathbf{P}[\xi\in P^\lambda]\bigr\}
\\
&\qquad\le\max\bigl\{\mathbf{P}[\xi\in P]-\mathbf{P}[\xi'\in P_{\lambda/2}],\, \mathbf{P}[\xi'\in P]-\mathbf{P}[\xi\in P_{\lambda/2}]\bigr\}
\\
&\qquad\le L_{N_m}\biggl(\mathcal{L}(\xi),\mathcal{L}(\xi'), \frac{\lambda}2\biggr).
\end{align*}
Since this holds for any $P\in\mathcal P_m$\ we arrive at inequality
$$
\pi_{m}(\,{\cdot}\,,{\cdot}\,, \lambda)\le L_{N_m}\biggl(\,{\cdot}\,,{\cdot}\,, \frac{\lambda}2\biggr).
$$

Thus, the second inequality of Theorem~\ref{th8} follows from the second inequality of Theorem~\ref{th2}. The constants depending on $N_m$ may be treated as constants depending on $m$. The first inequality follows from the second inequality by standard reasoning. Theorem~\ref{th8} is proved.\medskip

\noindent {\it Proof of Theorem\/ $\ref
 {th4}$.} Fix some  polyhedron $P\in\mathcal P_m$:
\begin{align*}
    P=\big\{x\in\mathbf H: \langle x,t_j\rangle\le b_j, \ j=1,\ldots, m\big\}.
\end{align*}where $t_j\in\mathbf H$, $\|t_j\|=1$, $b_j\in\mathbf R$,
$j=1,\ldots, m$. Let $\mathbf L_t\subset\mathbf H$ be the linear span of vectors
$$
\big\{t_j, j=1,\ldots, m\big\}, \quad k=\dim \mathbf L_t\le m,
$$
and let $\mathbb P_t:\mathbf H\to\mathbf L_t$ be the orthogonal projection operator on the subspace $\mathbf L_t$.
Consider the polyhedron $\overline{P}\subset\mathbf L_t$  defined as
\begin{align*}
 \overline   P=\big\{x\in\mathbf L_t: \langle x,t_j\rangle\le b_j, \ j=1,\ldots, m\big\}.
\end{align*}
It is easy to see that, for any random vector $\zeta\in\mathbf H$, we have
\begin{align*}
  \langle\mathbb   P_t\zeta,t_j\rangle=\langle \zeta,t_j\rangle,  \quad j=1,\ldots, m.
\end{align*}
Therefore,
 \begin{equation}
    \P[\zeta\in P]=  \P[\mathbb P_t\zeta\in \overline{P}]\quad\text{and}\quad\P[\zeta\in P_{\lambda}]=  \P[\mathbb P_t\zeta\in \overline{P}_{\lambda}],
\label{div7}
\end{equation}
where
$$
 \overline   P_{\lambda}=\big\{x\in\mathbf L_t: \langle x,t_j\rangle\le b_j+{\lambda}, \ j=1,\ldots, m\big\},\quad \lambda >0.
$$
Similarly, it is not difficult to show that   $ \{\zeta\in P^{\lambda}\}$ and $\{\mathbb P_t\zeta\in \overline{P}^{\lambda}\}$
are differemt descriptions of the same event. Therefore,
 \begin{equation}
    \P[\zeta\in P^{\lambda}]=  \P[\mathbb P_t\zeta\in \overline{P}^{\lambda}].
\label{div8}
\end{equation}
For a better understanding of the situation, it is useful to mentally consider the case where $ d = 3 $ and $ k = 2 $.

 The distributions of $k$-variate vectors $\mathbb P_t\xi,\mathbb P_t\eta\in\mathbf L_t$ actually satisfy the same $k$-dimensional conditions as the distributions of random vectors $\xi,\eta\in\mathbf H$ with compared infinite-dimensional distributions $F$ and $D^{**}$, $D^{*}$ or $D$ from Theorem \ref{th4}. In order to verify that, one should argue like in the proof of Theorem~\ref{th2} replacing operator $\mathbb A$ by operator $\mathbb P_t$ and using that $\|\mathbb P_t\|=1$.
 Applying Theorems~$\ref {th2}$ and~$\ref {th8}$, we obtain that, for any $\lambda >0$,
\begin{align}
&\max\bigl\{\mathbf{P}[\mathbb P_t\xi\in \overline{P}]-\mathbf{P}[\mathbb P_t\eta\in \overline{P}_{\lambda}],\, \mathbf{P}[\mathbb P_t\eta\in \overline{P}]-\mathbf{P}[\mathbb P_t\xi\in \overline{P}_{\lambda}]\bigr\}
\nonumber
\\
&\qquad \le c(m)\biggl( p+\exp \biggl( -\frac{\lambda}{c(m)\tau} \biggr) \biggr)
\label{div6}
\end{align}
and
\begin{align}
&\max\bigl\{\mathbf{P}[\mathbb P_t\xi\in \overline{P}]-\mathbf{P}[\mathbb P_t\eta\in \overline{P}^{\lambda}],\, \mathbf{P}[\mathbb P_t\eta\in \overline{P}]-\mathbf{P}[\mathbb P_t\xi\in \overline{P}^{\lambda}]\bigr\}
\nonumber
\\
&\qquad\le c(m)\biggl( p+\exp \biggl( -\frac{\lambda}{c(m)\tau} \biggr) \biggr).
\label{adiv6}
\end{align}
The statement of Theorem~\ref{th4} follows now from~\eqref{div7}--\eqref{adiv6}.  Theorem~\ref{th4} is proved.\medskip

In our results, we assume, for simplicity, that
$$ a_i=\int x\,U_{i}\{dx\}=0, \quad i=1,\dots ,n.$$
If we remove this assumption, then it will be valid again after replacing distributions $F_i$ by
distributions $F_iE_{-a_i}=(1-p_{i})U_{i}E_{-a_i}+p_{i}V_{i}E_{-a_i}$. Of course, $U_iE_{-a_i}$ is concentrated on the ball of larger radius $2\tau$, but this does not imply any change of the rate of infinitely divisible approximation if we are not interested in numerical values of constants. In particular, applying inequalities \eqref{div2}--\eqref{adiv},
we get the bounds
\begin{equation}
L_m(F,\overline D) \le
c(m)\,\big(p+\tau(|\ln \tau|+1)\big),
\label{div22}
\end{equation}
\begin{equation}
\pi_m(F,\overline D) \le
c(m)\,\big(p+\tau(|\ln \tau|+1)\big),
\label{adiv22}
\end{equation}
and
\begin{equation}
L_m(F,\overline D,\lambda) \le
c(m)\,\Big( p+\exp \Big( -\frac{\lambda
}{c(m)\,\tau }\Big) \Big) ,\quad\lambda>0,
\label{div11}
\end{equation}
\begin{equation}
\pi_m(F,\overline D,\lambda) \le
c(m)\,\Big( p+\exp \Big( -\frac{\lambda
}{c(m)\,\tau }\Big) \Big) ,\quad\lambda>0,
\label{adiv11}
\end{equation}
where
$$
{\overline D}=\prod_{i=1}^n E_{a_i}e(F_iE_{-a_i}).
$$
Clearly, it is easy to write the corresponding analogues of approximating distributions \eqref{d1} and \eqref{d2} with the same rate of approximation as in \eqref{div22}--\eqref{adiv11}.

The situation considered in Theorems~\ref{th2} and~\ref{th8} can be interpreted as a comparison of the sample containing
independent observations of rare events with the Poisson point process
which is obtained after a Poissonization of the initial sample (see~\cite{GZ18},~\cite{14}).

Indeed, let $Y_1, Y_2,\dots, Y_n$ be independent not
identically distributed elements of a measurable space $(\mathfrak
Y,\mathcal S)$.
Assume that the set $\mathfrak Y$ is represented as the union of two disjoint measurable sets: $\mathfrak Y= \mathfrak
Y_1\cup\mathfrak Y_2$, with \,$\mathfrak Y_1,\,\mathfrak Y_2\in \mathcal S$, $\mathfrak
Y_1\cap\mathfrak Y_2=\varnothing$. We say that the $i$-th  rare event occurs if $Y_i\in\mathfrak Y_2$.
Respectively, it does not occur if $Y_i\in\mathfrak Y_1$.

Let $f:\mathfrak Y\to\mathbf  R^d$ be a Borel mapping
and $F_i={\mathcal L}(f(Y_i))$, $i=1,2,\dots,n$. Then
 distributions $F_i\in\mathfrak F_d$ can be represented as mixtures \begin{equation}
F_i=(1-p_i)\,U_i+p_i\,V_i, \label{(11}\end{equation} where $U_i, V_i\in\mathfrak F_d$
are conditional distributions of vectors $f(Y_i)$ given
$Y_i\in\mathfrak Y_1$ and $Y_i\in\mathfrak Y_2$ respectively,
\begin{equation} 0\le p_i=\mathbf{P}\big\{Y_i\in\mathfrak Y_2\big\}
=1-\mathbf{P}\big\{Y_i\in\mathfrak Y_1\big\}\le1. \label{(12}\end{equation}
By definition, we deal with rare events whereas the quantity
\begin{equation}  p=\max_{1\le i\le
n}p_i\label{(13}\end{equation}
 is small. In other words, this is the case  if our rare events are sufficiently rare.

 Denote \begin{equation}
F=\prod_{i=1}^nF_i,\quad D=\prod_{i=1}^n e(F_i). \label{(15}
\end{equation}
The sum
\begin{equation} S=f(Y_1)+
\dots+ f(Y_n)\label{(1533}
\end{equation}
has the distribution $F$.
 It is easy to see that $D$ is the
distribution of \begin{equation}T=\sum_{i=1}^n\sum_{j=1}^{\nu_i}f(Y_{i,j}), \label{(135}\end{equation}
where $Y_{i,j}$ \tc{and $\nu_i$, $i=1,\dots,n$, $j=1, 2, \dots$, are  random elements in $\mathfrak Y$ which are
independent in aggregate such that
%and random variables respectively with
 ${\mathcal L}(Y_{i,j})={\mathcal L}(Y_{i})$ and
${\mathcal L}(\nu_i)=e(E_1)$.}
Clearly, $ e(E_1)$ is the Poisson distribution
with mean 1.

Thus, the sum $T$ is defined similarly to $S$, but
the initial sample ${\mathbf Y}=(Y_1, Y_2,\dots, Y_n)$ is replaced by
its Poissonized version ${\mathbf \Pi}=\bgl\{Y_{i,j}:$ $i=1,\dots,n$,
$j=1, 2, \dots,\nu_i\bgr\}$. Poissonization of the sample is known
as one of the most powerful tools in studying empirical processes.
The random set ${\mathbf \Pi}$ may be considered as a realization of the Poisson point
process on the space $\mathfrak Y$ with intensity measure
$\sum_{i=1}^n {\mathcal L}(Y_{i})$. The important property of the Poisson point process
is the space independence: for any pairwise disjoint sets $A_1,\ldots,A_m \in\mathcal S$,
the random sets $\mathbf \Pi\cap A_1,\ldots,\mathbf \Pi\cap A_m \subset\mathfrak Y$ are independent in aggregate.
\emph{}As a consequence, investigation of the Poisson point
process ${\mathbf \Pi}$ is much easier  than studying the sample $\mathbf  Y$.
One can use the independence property since  the theory of independent objects is much more elaborated.

 Let relations~$\eqref{(11}$--$\eqref{(135}$ be satisfied and let,
for some $\tau\ge0$,
\[
 U_{i}\big\{
\big\{ y\in \mathbf{R}^{d}:\left\| y\right\| \le \tau \big\}
\big\} =1, \quad i=1,2,\ldots,n,
\]
and the  $V_{i}\in\mathfrak F_d$ are arbitrary distributions.
Define
\begin{equation}
a_i=\int_{\mathbf  R^d} x\,U_i\{dx\},\quad i=1,2,\ldots, n.\label{(14}\end{equation} Denote
\begin{equation}T^*=\sum_{i=1}^n\Big(a_i+\sum_{j=1}^{\nu_i}\big(f(Y_{i,j})-a_i\big)\Big), \label{(177}\end{equation}
Then
\begin{equation}
{\overline D}=\mathcal{L}(T^*)=\prod_{i=1}^n E_{a_i}e(F_iE_{-a_i}), \label{(15f}\end{equation}
and \begin{equation}
T^*=T-\Delta,\quad\mbox{where } \Delta=\sum_{i=1}^n(\nu_i-1)\,a_i,
\label{(15fi}\end{equation}
 and $\nu_i$ are i.i.d.\  Poisson
with mean $1$.
\bigskip

Theorem~\ref{th2} implies the following assertions about the closeness of distributions $F$ and~$ { D}$, see \eqref{(15}.

\begin{theorem}\label{45}Let the above conditions   be
satisfied. Then, for any $m\in\mathbf{N}$, $\lambda>0$ and $P\in \mathcal P_m$  defined in \eqref{wq}, we have
\begin{align}
&\max\bigl\{F\{P\}-D\{P_{2\lambda}\},\, D\{P\}-F\{P_{2\lambda}\}\bigr\}
\nonumber
\\
&\qquad \le c(m)\biggl( p+\exp \biggl( -\frac{\lambda}{c(m)\tau }\biggr) \biggr) +\sum_{j=1}^{m}\mathbf{P}\{|\langle \Delta,t_j \rangle| \ge\lambda\},
\label{g955}
\end{align}
where the polyhedron $P_{2\lambda}$ is defined in \eqref{g9877}.
\end{theorem}

 \noindent {\it Proof}. Note that $P_{\lambda}\in \mathcal P_m$ and $(P_{\lambda})_{\lambda}=P_{2\lambda}$.
Using  \eqref{div11}, we see that
  \begin{equation}
 \max\Big\{F\{P\}-\overline D\{P_{\lambda}\},\;\overline D\{P_{\lambda}\}-F\{P_{2\lambda}\}\Big\}
\le c(m)\,\Big( p+\exp \Big( -\frac{\lambda
}{c(m)\,\tau }\Big) \Big). \label{0g955}\end{equation}By definition,
\begin{equation} \overline D\{P_{\lambda}\}=\mathbf P\big\{\langle T^*,t_j\rangle\le b_j+\lambda, \ j=1,\ldots, m\big\}, \label{60000}\end{equation}
\begin{equation} D\{P\}=\mathbf P\big\{\langle T,t_j\rangle\le b_j, \ j=1,\ldots, m\big\}, \label{900004}\end{equation}
\begin{equation} D\{P_{2\lambda}\}=\mathbf P\big\{\langle T,t_j\rangle\le b_j+2\lambda, \ j=1,\ldots, m\big\}\,. \label{800002}\end{equation}
Using \eqref{(15fi}, \eqref{60000}, \eqref{900004}, we obtain inequality \begin{equation} D\{P\}\le\overline D\{P_{\lambda}\} +\sum_{j=1}^{m}\mathbf{P}\big\{ \left|\langle \Delta,t_j \rangle \right| \ge\lambda \big\}.  \label{20000}\end{equation}Similarly, by  \eqref{(15fi}, \eqref{60000}, \eqref{800002}, we have\begin{equation} \overline D\{P_{\lambda}\}\le D\{P_{2\lambda}\}+\sum_{j=1}^{m}\mathbf{P}\big\{ \left|\langle \Delta,t_j \rangle \right| \ge\lambda \big\}.  \label{200006}\end{equation} Inequality \eqref{g955}\ follows now from \eqref{0g955}, \eqref{20000} and \eqref{200006}.Theorem \ref{45} is proved.

Theorem \ref{45} is a generalization of \cite[inequalities (15) and (16) of Theorem 9]{GZ18}.

The probabilities $\mathbf{P}\big\{ \left|\langle \Delta,t_j \rangle \right| \ge\lambda \big\}$ may be estimated using Bernstein's
inequality, see \cite[inequality (17)]{GZ18}.

\end{document}